\documentclass[11pt]{article}
\usepackage{setspace}
\usepackage{textcomp}
\usepackage{amsmath, amsthm, amssymb}
\usepackage[left=1.5in,top=1in,right=1in,nohead]{geometry}
\usepackage[dvips]{graphicx}
\usepackage{hyperref}
\usepackage{latexsym}
\usepackage{mathrsfs}
\usepackage{rawfonts}

\def\beg{\begin}

\newtheorem{thm}{Theorem}[section]

\newtheorem{rmk}[thm]{Remark}

\usepackage{makeidx}
\makeindex
\def\a{\alpha}

\def\bct{\begin{center}}
\def\ect{\end{center}}
\def\beg{\begin}

\def\<{\langle}
\def\>{\rangle}
\def\mbb{\mathbb}
\def\mbbc{\mathbb C}

\def\tn{\textnormal}

\newcounter{exmp}

\title{Geometric Invariant Theory and Einstein-Weyl Geometry}

\author{Mustafa Kalafat}

\begin{document}
\maketitle

\begin{abstract}
In this article, we give a survey of Geometric Invariant Theory for Toric
Varieties, and present an application to the Einstein-Weyl Geometry.
We compute the image of the Minitwistor space of
the Honda metrics as a categorical quotient according to the most efficient
linearization. The result is the complex weighted projective space
$\mathbb{CP}_{1,1,2}$. We also find and classify all possible quotients.
\end{abstract}


\section{Introduction} 

Let $(M,g)$ be a self-dual Riemannian $4$-manifold. This means that the
anti-self-dual Weyl tensor $W_-$ vanishes. In this case
\cite{ahs} construct a complex $3$-manifold $Z$ called the {\em
Twistor Space} of $M$, and 
a fibration by
holomorphically embedded rational curves. \vspace{.2cm}

\hspace{3.5cm} \begin{tabular}{cccc}
$\mathbb{CP}_1$ & $\to$ & $Z$ & Complex $3$-manifold \\
&&$\downarrow$& \\
&&$M^4$ & Riemannian $4$-manifold
\end{tabular}
\vspace{.2cm}

Suppose moreover that $M$ admits
a free isometric circle($S^1$) action. Then the quotient manifold
$M/S^1$ is naturally equipped with a so-called {\em Einstein-Weyl
Geometry}\index{Einstein-Weyl Geometry}. That is to say we have a
triple $(M/S^1,[h],D)$ where $[h]$ is a conformal class, here for
the induced metric of the quotient, and $D$ is a torsion-free affine
connection. The condition
\beg{equation}\tn{Ric}_{(ij)}=\lambda h_{ij}
\tag{Einstein-like} \end{equation}
more precisely $\tn{Ric}(u,v)+\tn{Ric}(v,u)=2\lambda h(u,v)$ and besides that the following
\beg{equation}Dh=\alpha\otimes h  
\tag{Weyl Connection} \end{equation} 
for some 1-form $\alpha$ are to be satisfied. 
This action can naturally be extended to a holomorphic $\mathbb C^*$-action over
the twistor space. We call the corresponding quotient $Z/\mathbb C^*$
 the $Minitwistor space$ \index{Minitwistor space} of the self-dual manifold.
 It is a very natural question to ask what is this 
quotient space.
We know that if the twistor space is algebraic or Moishezon the quotient 
becomes a complex surface with singularities in general.   


In the march of 2004, Honda gave an explicit description for
the twistor space of certain 
self-dual metrics on $3\mathbb{CP}_2$
admitting a free isometric circle action, equivalently a nowhere zero
Killing Field as follows.

\begin{thm}[Nobuhiro Honda,\cite{honda28}]\label{tw28}
Let $g$ be a  self-dual metric on $3\mathbb{CP}_2$  which admits a
non-trivial Killing Field. Suppose further that it is of positive scalar curvature type,
and not conformally equivalent to the hyperbolic ansatz self-dual metrics 
of LeBrun's  \cite{leexplicit}.\\ 
Then the twistor space is
a small resolution of 
the double cover of  $\mathbb{CP}_3$
branched along a quartic, equation of which is given in some homogeneous coordinates by
$$\left(Z_2 Z_3+Q(Z_0,Z_1)\right)^2-Z_0Z_1(Z_0+Z_1)(Z_0-aZ_1)=0$$
 where $Q(Z_0,Z_1)$ is a
quadratic form of $Z_0$ and
$Z_1$ with real coefficients, and $a\in \mathbb R^+$.\\
Moreover, the naturally induced real structure on $\mathbb{CP}_3$ is given by
$$\sigma(Z_0:Z_1:Z_2:Z_3)=\left(\bar{Z}_0:\bar{Z}_1:\bar{Z}_3:\bar{Z}_2\right),$$
and the naturally induced $U(1)$-action on
$\mathbb{CP}_3$ is given by
$$(Z_0:Z_1:Z_2:Z_3)\mapsto
\left(Z_0:Z_1:e^{i\theta}Z_2:e^{-i\theta}Z_3\right)~~for~~
e^{i\theta}\in U(1).$$ \end{thm}
\vspace{-.25cm}

To construct the minitwistor space of a Honda metric, we appeal
to the Geometric Invariant Theory (GIT) for Toric Varieties. This celebrated
theory was developed by D. Mumford around 1970's to understand the quotients of group
actions on manifolds. 
We compute the image under the double branched cover, so that we could
be able to recover the original minitwistor space by taking a double
cover along the related branch locus. GIT computes the quotients according to
some linearizations. It takes out some bad orbits, called the unstable orbits
and gives a toric variety as a result. 
We do computations for each linearization and finally figure the way to
minimize the number of unstable orbits. Summarizing our main Theorem \ref{quotient} and efficiency arguments in Section \ref{efficiencyandclassification} we have
obtained 

\vspace{1mm} 
\noindent {\bf Theorem A.} 
The image of the Minitwistor space of
a Honda metric in \cite{honda28} according to some efficient
specific linearization is the complex weighted projective
space $\mathbb{CP}_{1,1,2}$.

The idea is to compute the coordinate rings of the variety obtained,
and sketching the fan or the polytope of the toric variety to
realize an isomorphism with the fan or polytope of the
$\mathbb{CP}_{1,1,2}$. Yet, one can show that even this most refined
quotient excludes $\mathbb{CP}_1$-many orbits. But the GIT for Toric
Varieties does not provide a better solution than Theorem A. 
We define and discuss the efficiency and classification of quotients
arising from all possible linearizations. For our purposes, best linearizations
are the "efficient" ones as dicussed in Section \ref{efficiencyandclassification}. We are interested in the geometric(visual)
perspective, so that reducing the unstable orbits in terms of dimension,  measure or number of connected components is desirable for us. 
Summarizing the Theorems \ref{efficiencydimensions},  \ref{efficiencydimensions2} and Corollary \ref{efficiencyrestquotients},

\vspace{1mm} 
\noindent {\bf Theorem B.} The only possible categorical quotients of $\mathbb C^4$ under the $\mathbb C^{*2}$ action described by the matrix
$A=\left[\begin{array}{cccc}1&1&1&1\\ 0&0&1&-1 \end{array}\right]$ are the empty set, 
$\mathbb C$, $\mathbb{CP}_1$ and $\mathbb{CP}_{1,1,2}$.

Honda also describes these minitwistor spaces in \cite{hondamini28} but in a somewhat ad hoc way. So a GIT construction is desirable. 
The method we use can be applied to compute different
minitwistor spaces and also can be developed to be more effective. The project is to apply the general geometric invariant
theory and figuring out the complete information about these
quotients. 
It is a very common problem to figure out the
minitwistor spaces from the twistor spaces, and
there are a number of self-dual metrics for which the minitwistor
space is waiting to be computed. A systematic application of GIT
will address and solve many problems in the area. This paper should be considered as
a modest start for this program.
In \S\ref{secgit}-\S\ref{sector}
we give a review of GIT and Toric Varieties. 
Our survey 
owe much to 
the excellent resources \cite{dolgachevgit,brasselet} and \cite{mukai}.
Finally in \S\ref{secmini}-\S\ref{efficiencyandclassification} we
present our applications. \vspace{1mm}

\noindent {\bf Acknowledgments.} I want to thank to Claude LeBrun for his
directions, Alastair Craw for his lectures on the GIT and comments, Fr\'ed\'eric Rochon, Joel Robbin and the referee for careful examination, comments and corrections on the previous draft.

\section{Action of a torus on an affine space}\label{secgit}


In this section we will analyze the actions of the algebraic torus group $T=(\mbbc^*)^r$ on the affine space $\mbbc^n$ and understand the quotients
arisen this way. \vspace{1mm}

Recall that a {\em character} $\chi$ of an abelian group, with values in a field is
a homomorphism from the group to a multiplicative group, i.e.
satisfying $\chi(gh)=\chi(g)\chi(h)$. Moreover $\chi(T)$ stands for
the group of characters of $T$. We have the fact that any character
~$\chi : T \longrightarrow \mbbc^*$~
is given by \cite{dolgachevgit,mukai}
$$\chi(t)=\chi(t_1\cdots t_r)=t_1^{a_1}t_2^{a_2}\cdots t_r^{a_r}=\prod_{i=1}^r t_i^{a_i} $$
for $t_i \in \mbbc$, $a_i\in\mbb Z$. So we have the isomorphism
~$\chi(T)\approx \mbb Z^r$~ for the space of characters. \vspace{.05in}

Consequently, after diagonalization, a $T$ action on $\mbbc^n$ is
written as

$$t\cdot \left( \beg{array}{@{}c@{}}Z_1\\ \vdots \\ Z_n \end{array}
\right)=\left( \beg{array}{@{}c@{}}\chi_1(t)Z_1\\ \vdots \\ \chi_n(t)Z_n
\end{array} \right)=\left( \beg{array}{@{}c@{}}t^{a_1}Z_1\\ \vdots \\ t^{a_n}Z_n \end{array}
\right)=\left( \beg{array}{@{}c@{}}t_1^{a_{11}}\cdots t_r^{a_{r1}}Z_1\\ \vdots \\
t_1^{a_{1n}}\cdots t_r^{a_{rn}}Z_n
\end{array} \right),\vspace{2mm}$$
so the matrix $A=[a_{ij}]\in M_{r\times n}(\mbb Z)$ encodes the
action. \vspace{2mm}


More generally, let $\sigma : T\times X \to X$ be an action of the group $T$ on the complex manifold $X$ by complex automorphisms. For a holomorphic line
bundle $\pi : L\to X$,  we define

\beg{defn} A {\em linearization}\index{linearization} of the holomorphic line bundle $L$ with respect to the action of $T$ is an action $\overline{\sigma} : T\times L \to L$ so that \hspace{1mm} \vspace{1mm}\\
($1$) The following diagram commutes \\
$$\beg{array}{ccc} T\times L & \stackrel{\overline{\sigma}}{\longrightarrow} &
L\\ \llap{ {\tiny $id\times\pi$}  } {\huge \downarrow} &&  \downarrow \rlap{{\tiny $\pi$}} \\
T\times X & \stackrel{\sigma}{\longrightarrow} & X
\end{array}$$
($2$) The zero section $X\approx L_0  \subset L$ is $T$- invariant.
\end{defn}
So this is the extension of the action $\sigma$ to $L$, preserving
the fibers, i.e. points on a fiber map  onto the same fiber under
the action of an element. It follows from the definition that this
action on a fiber $\overline{\sigma}_t : L_p \to L_{tp}$ for any
$t\in T$ and any $p\in X$ is a linear isomorphism.               

In our case, the action of $\mbbc^{r*}$ on $\mbbc^n$ is given by the
matrix $A=(a_1\cdots a_r)\in M_{n\times r}(\mbb Z)$. Consider the
trivial line bundle $\underline{\mbbc}\to \mbbc^n$. Fix
$\alpha=(\alpha_1\cdots\alpha_r)\in \mbb Z^r$. Extend the action
over to the bundle $\underline{\mbbc}$ as follows
$$t\cdot(Z,W)=(t\cdot Z,t^\alpha W)=(t\cdot Z,t_1^{\alpha_1}t_2^{\alpha_2}\cdots t_r^{\alpha_r}W)
~~\textnormal{where}~~Z\in\mbbc^n , W\in \mbbc.$$ We denote this
linearized line bundle by $L_\alpha$. So any $a\in \mbb Z^r$ gives
an extension or a linearization.

Recall that the holomorphic sections of the trivial line bundle are identified
with the polynomials $F\in\mbbc[Z_1\cdots Z_n]$, like the homogenous
polynomials for bundles over $\mbb P^n$. A section $F$ is an
invariant section of $L_\alpha$ if $$t\cdot(Z,F(Z))=(t\cdot
Z,t^\alpha\cdot F(Z))=(t\cdot Z,F(t\cdot Z)),$$
which amounts to
$$t^\alpha\cdot F(Z)=F(t\cdot Z),$$ that is $$t_1^{\alpha_1}\cdots
t_r^{\alpha_r}F(Z_1\cdots Z_r)=F(t^{a_1}Z_1\cdots t^{a_r}Z_r).$$ 
The action of $\overline{\sigma}$ on $L$ induces an action on
$L^{\otimes d}$ as for a decomposable $ l\in L_p^{\otimes d}$,
$$\overline{\sigma}_t(l)=\overline{\sigma}_t(l_1\otimes\cdots\otimes
l_d)=\overline{\sigma}_t(l_1)\otimes\cdots\otimes
\overline{\sigma}_t(l_d) \in L_{tp}^{\otimes d}.$$
Likewise, $G$ is an invariant section of $L_\a^{\otimes d}$ if for
$G=F_1\cdots F_d$
$$\beg{array}{rcl} G(t\cdot Z)&=&F_1(t\cdot Z)\cdots
F_d(t\cdot Z)\\ &=&(t^\alpha\cdot F_1)\cdots (t^\alpha\cdot F_d)\\
&=&t^{\alpha d}\cdot F_1\cdots F_d\\ 
&=&t^{\alpha d}\cdot G(Z).\end{array}$$
Imposing the above condition one proves that

\beg{prop}[\cite{dolgachevgit}] $G \in H^0(\mbbc^n, L^{\otimes d}_\a)^T$ i.e. $G$ is an
invariant section of the linearized line bundle $L_\a^{\otimes d}$
~iff~ it is a linear combination of monomials $Z^m=Z_1^{m_1}\cdots
Z_n^{m_n}$  such that
\begin{equation}\label{monomialequation} \tag{Monomial Equation}
[A ~,-\alpha ~] \left[\beg{array}{@{}c@{}}m\\
d\end{array}\right] = 0_r \end{equation}
where $A \in M_{r\times n}(\mbb Z)$
is the action matrix, $\alpha\in \mbb Z^r$ is the tuple for the
extension.
\end{prop}

\beg{proof} Say $G=Z^m$ , then \beg{eqnarray} G(t\cdot
Z)&=&t^{\alpha d}\cdot G(Z)\nonumber\\
G(t^{a_1}Z_1\cdots t^{a_n}Z_n)
&=&(t_1^{\alpha_1}\cdots t_n^{\alpha_n})^dZ^m\nonumber\\
(t^{a_1}Z_1)^{m_1}\cdots (t^{a_n}Z_n)^{m_n}&=& ~~~~~''\nonumber\\
t^{a_1 m_1}\cdots t^{a_n m_n}Z^m&=& ~~~~~''\nonumber\\
(t_1^{a_{11}}\cdot\cdot t_r^{a_{r1}})^{m_1}\cdots
(t_1^{a_{1n}}\cdot\cdot t_r^{a_{rn}})^{m_n} Z^m&=& ~~~~~'' \nonumber
\end{eqnarray}
Comparing the powers of $t_i$'s from both sides we obtain the
equality
\beg{eqnarray*} a_{i1}m_1+\cdots +a_{in}m_n &=&\alpha_i d, \\
\left[ a_{i1} \cdots a_{in} \right]\left[\beg{array}{@{}c@{}}m_1\\
\vdots \\ m_n
\end{array}\right]&=&[\alpha_i]\, d ~~~\textnormal{for any}~~ 1\leq i\leq r,\\
Am&=&\alpha d. 
\end{eqnarray*}
\end{proof}

\beg{ex} Consider the following action  of $\mbbc^{*2}$ on
$\mbbc^4$, 
 
$$(t_1,t_2)\cdot \left( \beg{array}{@{}c@{}}X\\ Y\\ Z\\ W
\end{array} \right)=\left(
\beg{array}{@{}c@{}}t_1X\\ t_1^{-n}t_2Y\\ t_1Z\\ t_2W
\end{array} \right)~~,~~\alpha =\left(\beg{array}{@{}c@{}}1\\ 1\end{array} \right).$$
The action matrix is $A=\left[ \beg{array}{cccc}1&-n&1&0\\
0&1&0&1\end{array} \right]$ and the monomials for the invariant
sections are obtained from the equation
$\left[ \beg{array}{cccc|c}1&-n&1&0&-1\\
0&1&0&1&-1\end{array} \right]\left[ \beg{array}{@{}c@{}}m_1\\ m_2\\ m_3\\
m_4\\ d
\end{array} \right]=0_2.$
\end{ex}


Next we are going to give some definitions in the Geometric
Invariant Theory (GIT), which deals with the actions of groups on
manifolds, and figuring out their corresponding quotients.

\beg{defn}[Stability\cite{dolgachevgit}] Let $L$ be a
$T$-linearized line bundle on the
algebraic variety $X$ and let $x\in X$, then\\
(i)~$x$ is called {\em semi-stable}\index{semi-stable} with respect
to $L$ if it belongs to the set $
X\backslash\{s=0\}\subset
\mbbc^n $ (affine) for some $m>0$ and some $s\in H^0(X,L^m)^T$.\\
(ii)~$x$ is called {\em unstable} \index{unstable} with respect to
$L$ if it is not
semi-stable. 
\end{defn} 
We respectively denote by $X^{ss}(L)$ and $X^{us}(L)$, the set of
semi-stable and unstable points in $X$.

\beg{defn}[Categorical Quotient\cite{dolgachevgit}]
\index{categorical quotient} A {\em categorical quotient} of a
$T$-variety $X$ is a $T$-invariant morphism $p : X\to Y$ such that
for any $T$-invariant morphism $g : X\to Z$, there exist a unique
morphism $\bar{g} : Y\to Z$ satisfying $\bar{g} \circ p = g$. $Y$ is
written sometimes as $X/\!\!/T$ and also called the {\em categorical
quotient}. 
 \end{defn}

The GIT guarantees a (good) categorical quotient
$X^{ss}(L_\alpha)/T$, see \cite{dolgachevgit} (pp.118), denoted
alternatively by $X(L)/\!\!/\!_\alpha T$. This is the quotient
obtained by taking out the unstable orbits. So according to the GIT,
semi-stable points has this well behaving quotient described as
follows.

\beg{prop}[\cite{dolgachevgit}] If $X$ is projective and $L$ is
ample, we can compute the categorical quotient by
$$X(L)/\!\!/\!_\alpha T=\tn{Proj}\left( \bigoplus_{d\geq 0}H^0
(X,L_\alpha^{\otimes d})^T \right).$$
\end{prop}

\section{Toric Varieties}\label{lebrunkim}\label{sector}
Let $V\subset \mbb
C^n$ be an affine variety. We define its  {\em (affine) Coordinate
ring}\index{coordinate ring} to be
$$\mbbc[V]=\mbbc[z_1 \cdots z_n]|_V.$$
This is to say the coordinate ring is the ring of {\em regular
functions} according to the terminology of \cite{shafarevich}. If
we look at the restriction map
$$restr : \mbbc[z_1 \cdots z_n]\longrightarrow \mbbc[z_1 \cdots
z_n]|_V$$ we see that its kernel is equal to $I_V$, the vanishing
ideal of $V$. So the coordinate ring becomes $$\mbbc[V]=\mbbc[z_1
\cdots z_n]/I_V.$$
For any ring $R$, we define its {\em maximal spectrum}\index{maximal spectrum} by
$$\tn{Specm}(R)=\{I < R : I~\textnormal{is a maximal ideal}~ \}.$$
For any affine variety $V\subset\mbb C^n$, defining the Zariski Topology on each side we have the homeomorphism $V\approx \tn{Specm}(\mbbc[V])$
between an affine variety and the maximal spectrum of its coordinate ring. As the trivial case,
$\mbbc^n\approx \tn{Specm}\,\mbb C[z_1 \cdots z_n]$, where a point $a\in \mbbc^n$ corresponds to its vanishing ideal $$I_{\{a\}}=\mbbc[z](z_1-a_1)+\cdots +\mbbc[z](z_n-a_n)=\<z_1-a_1,\cdots,z_n-a_n\>.$$ The maximal ideals of the latter type consumes the maximal ideals of the polynomial ring $\mbbc[z_1\cdots z_n]$, see \cite{mukai}, which is referred
as the Weak Nullstellensatz in the literature \cite{brasselet}.
The full spectrum is the larger space of prime ideals with which we do not deal here.\vspace{1mm}

For any group $G$, the {\em group ring} $\mbbc[G]$ is the vector space with basis
$\{[g]\}_{g\in G}$ together with the bilinear product based on group multiplication. This amounts to a {\em $\mbbc$-algebra}. If we relax the inverse condition on a group then we get similar operations and obtain the {\em monoid algebra}. As an example we can consider $\mbbc[\mathbb Z^n]$ which is the same 
as the algebra of Laurent polynomials $\mbbc[Z_1^{\pm1}\cdots Z_n^{\pm n}]$ under the correspondence $m\in \mathbb Z^n$ to $Z^m=Z_1^{m_1}\cdots Z_n^{m_n}$. Similarly if we have a submonoid of $\mathbb Z^n$, its monoid algebra will be the subalgebra of Laurent polynomials generated by the corresponding monomials. Notice that we are using the same notation with the coordinate ring. The reader 
is expected to 
interpret the meaning from the context, will depend on what lies in the bracket, a geometric object or an algebraic one.\vspace{1mm}

We first go into the definition of an affine toric variety. For that
purpose we take a cone $\sigma$ in $\mbb R^n$ satisfying the
conditions of the following definition for the canonical lattice
$N\approx\mbb Z^n\subset \mbb R^n.$
\beg{defn}
[Cone] Let $A=\{x_1\cdots
x_r\}\subset\mbb R^n$ be a finite set of vectors. Then
\beg{itemize}

\item  The set $\sigma=\{x\in\mbb R^n :
x=\lambda_1x_1+\cdots +\lambda_rx_r~,~\lambda_i\geq0\}$
 is called a \emph{
 cone}.

\item $\sigma$ is called a \emph{lattice cone} if all the vectors $x_i\in A$ belong to $N$.

\item $\sigma$ is called \emph{strongly convex} if it does not contain any straight line going through the origin, i.e. $\sigma\cap -\sigma=\{0\}$.
\end{itemize}\end{defn}
\noindent In this case we define the {\em affine toric variety} corresponding to $\sigma$
as $$U_\sigma:=\tn{Specm}\,\mbb C[\check{\sigma}\cap N^*]$$
where the dual is defined to be $\check{\sigma}=\{ u\in\mbb R^n : \<u,\sigma\>\geq 0  \}$ and $N^*=Hom_{\mbb Z}(N,\mbb Z)$. By abuse of notation, one can also write $U_\sigma=\tn{Specm}\,\mbbc[\check{\sigma}]$. 
Similar to the way that the cones correspond to an affine toric variety,
some collection of cones called fans correspond to a toric variety.
More precisely
\beg{defn}[Fan,\cite{brasselet}] A \emph{fan} $\Delta$ is a finite union of cones such that \beg{itemize}
\item The cones are 
lattice and strongly convex.
\item Every face of a cone of $\Delta$ is again a cone of $ \Delta$.
\item  $\sigma\cap\sigma '$ is a common face of the cones $\sigma$ and $\sigma '$ in $\Delta$.
\end{itemize}\end{defn}

Now for a fan $\Delta$ in $N$, we can naturally glue $\{U_\sigma: \sigma\in\Delta\}$ together to obtain a Hausdorff complex analytic space
$$X_\Delta :=\bigcup_{\sigma\in\Delta}U_\sigma$$
which is irreducible and normal with dimension equal to $\tn{rank}(N)$ and called
the {\em Toric Variety}\index{toric variety} \cite{oda} associated to the fan $(N,\Delta)$. It is topologically endowed with an open cover by the affine toric varieties $U_\sigma=\tn{Specm}\, \mbbc[\check \sigma]$.

Summarizing what we did in high brow terms \cite{dolgachevgit}, we constructed the $U_\sigma=\tn{Specm}\,\mbb C[\check{\sigma}\cap N^*]$
as the affine variety with $\mbbc[U_\sigma]$ isomorphic to
$\mbb C[\check{\sigma}\cap N^*]$. Since for any $\sigma,\sigma
'\in\Delta$, $\sigma\cap\sigma '$ is a face in both cones, we obtain
that 
$\mbb C[(\sigma\cap\sigma ')\check{} \cap N^*]$
is a localization
of each algebra $\mbb C[\check{\sigma}\cap N^*]$ and $\mbb C[\check{\sigma}'\cap N^*]$. This shows that $\tn{Specm}\,\mbb C[(\sigma\cap\sigma
')\check{}\cap N^*]$ is isomorphic to an open subset of $U_\sigma$ and
$U_{\sigma '}$. Which allows us to glue together the varieties
$U_\sigma$'s to obtain the toric variety $X_\Delta$.\vspace{1mm}

Returning to our case where we have an action of a torus $T$ on an affine space $\mbbc^n$, $\alpha$-linearized over to a line bundle $L$, we will now produce a fan and a toric variety out of this linearization. Notice that we have a 
natural isomorphism of graded algebras
$$\bigoplus_{d\geq 0}H^0(\mbbc^n,L_\alpha^{\otimes d})^T\approx \mbbc[S]=\bigoplus_{d\geq 0}\mbbc[S_d]$$
where $S$ is the monoid of elements $m \in \mathbb Z^n$ solving the 
\ref{monomialequation}. $\mbbc[S_d]$ is the linear span of $S_d$
which is the set of d-th solutions of the \ref{monomialequation}. 
It is degree-d homogenous part of the finitely generated $\mbbc[S]$. 
The ideal $$\mbbc[S]_{>0}:=\bigoplus_{d>0}\mbbc[S_d]=\<Z^{m_1},\cdots, Z^{m_s}\>$$ is finitely generated by a minimal set of monomial generators where 
$m_j=(m_{1j}\cdots m_{nj})$. For $I_j:=\{i~|~m_{ij}\neq 0\}$ and $Z_I:=\Pi_{i\in I}Z_i$ where $I\subset\{1\cdots n\}$ we have the equality
$$D(Z^{m_j}):=\mbbc^n-\{Z^{m_j}=0\}=\mbbc^n-\{Z_{I_j}=0\}=:D(Z_{I_j}).$$
By its definition the semi-stable locus becomes $$(\mbbc^n)^{ss}(L_\alpha)=\bigcup_{i=1}^s D(Z_{I_j}).$$
Thinking the matrix as a map $A:\mbb Z^n\longrightarrow \mbb Z^r$, let $M=\tn{Ker}A\subset \mathbb Z^n$ and for $1\leq j \leq s$ define
$$R_j:=\mbbc[D(Z_{I_j})]^T=\left\{ {F(Z)\over Z_{I_j}^p}~:~p\geq 0 \tn{ and } 
F(Z)\in Z_{I_j}^p\mbbc[M] \right\}.$$
Next, we will be gluing together some affine varieties coordinate rings 
of which are $R_j$'s. $M$ is a free abelian group of rank $l=n-\tn{rank}A$. 
Consider the map $(\mbb Z^n)^*\longrightarrow N=M^*$, which is given by restricting the linear functionals. Let $\{e_i\}$ be a basis for $\mbb Z^n$, 
$\{e_i^*\}$ be the dual basis with respect to the Euclidean metric, and 
$\{\overline e_i^*\}$ be their image in $M^*$. We define the convex cones $\sigma_{I_j}$'s or more concisely $\sigma_j$'s as the following span :
$$\sigma_j:=\<\,\overline e_i^*~|~i\notin I_j\,\>\subset N_{\mbb R}:=
N\otimes\mbb R \approx \mbb R^l.$$
One can show that $R_j\approx \mbbc[\check\sigma_j\cap M]$. $\sigma_j$'s form a fan $\Delta$, and this fan gives the toric variety we are seeking as the quotient, see \cite{dolgachevgit} for details. Consequently we have the folowing,

\beg{thm}[\cite{dolgachevgit}] Let $(\mbb Z^n)^*\longrightarrow M^*$ be the transpose of the inclusion $M\hookrightarrow \mbb Z^n$ and $N$ be its image. 
Let $\Delta$ be the $N$-fan
formed by the cones $\sigma_j , j=1\cdots s$ defined as above. Then
$$\mbb C^n(L)/\!\!/\!_\alpha T=(\mbb C^n)^{ss}(L_\alpha)/ T\approx X_\Delta.$$

\end{thm}

\beg{ex}\label{wps112} The weighted projective space $\mbb{CP}_{1,1,2}$ is by definition the quotient of $\mbb C^3-0$ by the
$\mbb C^*$-action given by the matrix $A=[1,1,2]$.

If we linearize the trivial bundle over $\mbb C^3$ by $\alpha=2$,
the linear system $Am=\alpha$ is just $a+b+2c=2$, and nonnegative
solutions for the triple $(a,b,c)$ are generated by
$$(2,0,0)~~ (1,1,0)~~ (0,2,0)~~ (0,0,1)$$
so that the coordinate rings are obtained as\vspace{4mm}

{\large
$\mbb C[\mbb N^4\cap \pi^{-1}(2)]=\mbb C[X^2,XY,Y^2,Z]$\vspace{3mm}

$\mbb C[U_1/\mbb C^*]=\mbb C[1,{Y\over X},{Y^2\over X^2},{Z\over
X^2}]=\mbb C[{Y\over X},{Z\over X^2}]=\mbbc[a,b]$\vspace{3mm}

$\mbb C[U_2/\mbb C^*]=\mbb C[{X\over Y},1,{Y\over X},{Z\over
XY}]=\mbb C[{X\over Y},{Y\over X},{Z\over
XY}]=\mbbc[a^{-1},a,a^{-1}b]$\vspace{3mm}

$\mbb C[U_3/\mbb C^*]=\mbb C[{X^2\over Y^2},{X\over Y},1,{Z\over
Y^2}]=\mbb C[{X\over Y},{Z\over
Y^2}]=\mbbc[a^{-1},ba^{-2}]$\vspace{3mm}

$\mbb C[U_4/\mbb C^*]=\mbb C[{X^2\over Z},{XY\over Z},{Y^2\over
Z},1]=\mbb C[{X^2\over Z},{XY\over Z},{Y^2\over
Z}]=\mbbc[b^{-1},ab^{-1},a^2b^{-1}]$}\vspace{3mm}\hspace{1mm}\\ if we
assign $a={Y\over X}$ and $b={Z\over X^2}$.\vspace{2mm}\\
Then since 
{\large $$\begin{array}{rcl}  \bigcup_{i=1}^4U_i & = & \mbb C^3 - \{ \{X^2=0\}
\cap \{XY=0\} \cap \{Y^2=0\} \cap \{Z=0\} \}   \\
&=& \mbb C^3 - \{ \{X=0\} \cap \{XY=0\} \cap \{Y=0\} \cap \{Z=0\}
\} \\
&=&\mbb C^3 - \{X=Y=Z=0\}  \end{array}$$ }
these are the coordinate
rings of the stated weighted projective space.\\
The moment polytope looks like :\newline

\begin{center}

\setlength{\unitlength}{2.5cm}
\begin{picture}(2,1)
\put(0,0){\line(1,0) {2} }%
\put(2,0){\line(-2,1){2} }%
\put(0,1){\line(0,-1) {1}} %

\put(0,0){\circle*{.05}} \put(-.25,-.25){$X^2$}

\put(1,0){\circle*{.05}} \put(.90,-.25){$XY$}

\put(2,0){\circle*{.05}} \put(2,-.25){$Y^2$}

\put(0,1){\circle*{.05}} \put(-.25,1){$Z$}

\end{picture}

\end{center}

\end{ex}\vspace{1cm}

\section{Minitwistor Space}\label{secmini}

The image of the Honda Minitwistor space (\ref{tw28}) is the
quotient of $\mbb{CP}_3$ by the $\mbb C^*$ action
$$(Z_0:Z_1:Z_2:Z_3)\mapsto
\left(Z_0:Z_1:\lambda Z_2:\lambda^{-1}Z_3\right)~~\tn{for}~~ \lambda\in
\mbb C^*.$$ On the other hand, to obtain $\mbb{CP}_3$, we already
have the classical $\mbb C^*$ action
$$(Z_0:Z_1:Z_2:Z_3)\mapsto
\left(\lambda Z_0:\lambda Z_1:\lambda Z_2:\lambda Z_3\right)~~\tn{for}~~
\lambda\in \mbb C^*.$$ Combining the two, the image equals to the
quotient of the $\mbb C^{*2}$ action by the matrix
$$A=\left[\beg{array}{cccc}1&1&1&1\\0&0&1&-1\end{array}\right]$$
on $\mbb C^4$. Now, extend this action to the trivial
line bundle over $\mbb C^4$. Choices are the linearizations. 
Among all of them, one of 
has the minimal number of unstable orbits.

 \beg{thm}\label{quotient} The categorical quotient of $\mbb C^4$ under the $\mbb
C^{*2}$ action described by the matrix
$$A=\left[\beg{array}{cccc}1&1&1&1\\0&0&1&-1\end{array}\right]$$
linearized by $\alpha=(2,0)$ is the weighted projective space $\mbb{CP}_{1,1,2}$.  \end{thm}

\beg{proof} The linear system $A m=\alpha$ is
$$\left. \beg{array}{rcl}a+b+c+d&=&2\\c-d&=&0\end{array}\right\}~~\tn{or}~~
 \left\{ \beg{array}{rcl}a+b+2d&=&2\\c&=&d\end{array}\right.$$
looking for nonnegative solutions, $1,0$ are the only
possibilities for $d$ since from the first equation $2d\leq 2$. So \begin{itemize}
\item $d=0 : a+b=2 , ~c=0 ~~\textnormal{yields the solutions}~~
(2~~0~~0~~0), (1~~1~~0~~0), (0~~2~~0~~0)$.
\item $d=1 : a+b=0 , ~c=1 ~~\textnormal{yields the solution}~~~
(0~~0~~1~~1)$. 
\end{itemize}
So the coordinate rings are\vspace{4mm}

{\large 
$\mbb C[\mbb N^4 \cap \pi^{-1}(2,0)]=\mbb
C[X^2,XY,Y^2,ZW]$\vspace{3mm}

$\mbb C[U_1/\mbb C^{*2}]=\mbb C[1,{XY \over X^2},{Y^2 \over
X^2},{ZW \over X^2}]=\mbb C[{Y \over X},{Y^2 \over X^2},{ZW
\over X^2}]=\mbb C[{Y \over X},{ZW \over X^2}]$\vspace{3mm}

$\mbb C[U_2/\mbb C^{*2}]=\mbb C[{X^2 \over XY},1,{Y^2 \over
XY},{ZW \over XY}]=\mbb C[{X \over Y},{Y \over X},{ZW \over
XY}]$\vspace{3mm}

$\mbb C[U_3/\mbb C^{*2}]=\mbb C[{X^2 \over Y^2},{XY \over
Y^2},1,{ZW \over Y^2}]=\mbb C[{X^2 \over Y^2},{X \over Y},{ZW
\over Y^2}]=\mbb C[{X \over Y},{ZW \over Y^2}]$\vspace{3mm}

$\mbb C[U_4/\mbb C^{*2}]=\mbb C[{X^2 \over ZW},{XY \over
ZW},{Y^2 \over ZW},1]=\mbb C[{X^2 \over ZW},{XY \over ZW},{Y^2
\over ZW}]$\vspace{4mm}\hspace{1mm}\\ } 
and these coordinate rings are
isomorphic to the ones for the $\mbb{CP}_{1,1,2}$ as in
(\ref{wps112}).
Realize the isomorphism by assigning $c={Y\over X} ,
d={ZW\over X^2}$ so that the coordinate rings respectively becomes
$$\mbbc[c,d]~,~\mbbc[c,c^{-1},c^{-1}d]~,~\mbbc[c^{-1},c^{-2}d]~,~
\mbbc[d^{-1},cd^{-1},c^2d^{-1}].$$ Besides, the moment polytope may
help to visualize this isomorphism :\newline\newline

\begin{center}

\setlength{\unitlength}{2.5cm}
\begin{picture}(2,1)
\put(0,0){\line(1,0) {2} }%
\put(2,0){\line(-2,1){2} }%
\put(0,1){\line(0,-1) {1}} %

\put(0,0){\circle*{.05}} \put(-.25,-.25){$X^2$}

\put(1,0){\circle*{.05}} \put(.85,-.25){$XY$}

\put(2,0){\circle*{.05}} \put(2,-.25){$Y^2$}

\put(0,1){\circle*{.05}} \put(-.25,1.1){$ZW$}
\end{picture}

\end{center}\vspace{.5cm}
\end{proof}
\hspace{-6mm}Realize that the union of $U_i$'s does not cover $\mbb C^4$ since
{\large
$$\begin{array}{rcl}\bigcup_{i=1}^4U_i &=& \mbb
C^4 - \{ \{X^2=0\} \cap \{XY=0\} \cap \{Y^2=0\} \cap
\{ZW=0\} \} \\
&=&\mbb C^4 - \{ \{X=Y=Z=0\} \cup \{X=Y=W=0\} \}. \end{array}$$}
\noindent Consequently, the points $[0:0:0:1],[0:0:1:0]$ in $\mbb {CP}_3$ are omitted in this quotient.



\section{Efficiency and Classification}\label{efficiencyandclassification}

In this section we analyze the efficiency of the linearization in
Theorem \ref{quotient}. Our notion of efficiency is based on the
maximum dimension of the omitted part under the action which we call the {\em efficiency dimension}. If this dimension is smaller, we say that the
corresponding linearization is more efficient. If two linearizations have the same efficiency dimension, then we consider the measure or number of connected components of the omitted piece to decide which one is more efficient. 
As an example, the linearization in the theorem 
has efficiency dimension $\tn{Ed}(2,0)=1$.

\beg{thm}\label{efficiencydimensions} Let $\alpha=(x,y) \in \mbb N^2$, 
i.e. a linearization. Then the following holds.
\beg{itemize}
\item If $y=0$    then the efficiency dimension $\tn{Ed}(0,0)=4$, $\tn{Ed}(1,0)=2$,~ moreover 
$\tn{Ed}(2m,0)=1$ and $\tn{Ed}(2m+1,0)=2$ for $m\geq1$.
\item If $y\geq1$ then the efficiency dimension $\tn{Ed}(x,y)\geq 3$.
\end{itemize}
\end{thm}

\beg{proof} Recall that we are considering the following system.
$$\left. \beg{array}{rcl}a+b+c+d&=&x\\c-d&=&y\end{array}\right\}~~\tn{or}~~
 \left\{ \beg{array}{rcl}a+b+2d&=&x-y\\c&=&d+y\end{array}\right.$$
Since we are concerned with nonnegative solutions, we need to have $x\geq y$ for that purpose.
Suppose first that $y=0$.
\beg{itemize}
\item $x=0$ : Solution space $SS=\{(0~0~0~0)\}$ is trivial. Charts are empty and
$\tn{Ed}=4$.

\item $x=1$ : $d=0$, $a+b=1,$ $c=0$. $SS=\<(1~0~0~0),(0~1~0~0)\>_+.$ 
The quotient turns out to be a $\mbb{CP}_1$ for the 
coordinate rings are

{\large
$\mbb C[\mbb N^4 \cap \pi^{-1}(1,0)]=\mbb
C[X,Y]$\vspace{3mm}

$\mbb C[U_1/\mbb C^{*2}]=\mbb C[{X\over Y}]=\mbb C[\beta]$ \vspace{3mm}

$\mbb C[U_2/\mbb C^{*2}]=\mbb C[{Y \over X}]=\mbb C[\beta^{-1}]$. }

The omitted locus $\{X=Y=0\}\subset \mbb{C}^4$ has dimension $\tn{Ed}=2$.

\item $x=2m$, ~$m\geq 1$ : We have $a+b=2(m-d)$ and $c=d$ in this case.
$a+b$ decreases evenly as $d$ increases. So the coordinate ring is as follows.
$$\mbb C[X^{2m},X^{2m-1}Y \cdots XY^{2m-1},Y^{2m}, 
   \{X^{2(m-1)}, X^{2(m-1)-1}Y \cdots Y^{2(m-1)}\}ZW  
\cdots ZW^m].$$
This suggests that the ommited locus $X=Y=Z=0$ or $X=Y=W=0$, which implies
that $\tn{Ed}=1$ for this case.

\item $x=2m+1$, $m\geq 1$ : We have $a+b=2(m-d)+1$ and $c=d$ in this case.
$d\leq m$ for a positive solution to exist.

The coordinate ring

$$\mbb C[X^{2m+1},X^{2m}Y \cdots XY^{2m},Y^{2m+1}, 
   \{X^{2m-1}, X^{2m-2}Y \cdots Y^{2m-1}\}ZW,$$
$$\{X^{2m-3}\cdots Y^{2m-3}\}ZW^2  \cdots \{X, Y\}ZW^m]$$

yields the omitted locus $X=Y=0$ hence the $\tn{Ed}=2$ in this case.
\end{itemize}
Now suppose $y\geq1$. Then from the second equation we have $c=d+y\geq1$. 
This tells us that $c$ is nonzero, consequently the hyperplane $Z=0$ always lies in the ommited locus.
\end{proof}

\beg{thm}\label{efficiencydimensions2} Let $\alpha=(x,y) \in \mbb N^2$, i.e. a linearization. Let $y\geq 1$.
Then we have the following dimensions and quotients. \beg{itemize}
\item If $x-y<0$ then $\tn{Ed}=4$, and the quotient is empty.   
\item If $x-y=0$ then $\tn{Ed}=3$, and the quotient is $\mbbc$. 
\item If $x-y=1$ then $\tn{Ed}=3$, and the quotient is a complex projective line $\mbb{CP}_1$.
\item If $x-y\geq 2$ then $\tn{Ed}=3$, and the quotient is the weighted
projective space $\mbb{CP}_{1,1,2}$.
\end{itemize} \end{thm}

\begin{figure}[!h]            
\bct  \includegraphics[width=\textwidth]{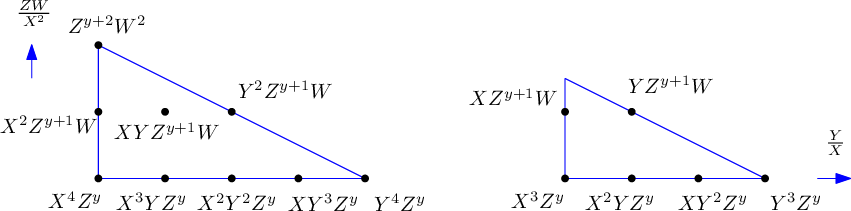}\\ \ect
\caption{\em\small Typical even$(x-y=4)$ and odd$(x-y=3)$ cases.} 
\label{ygeq1} \end{figure}

\beg{proof} The first three cases are similar to that of in the proof of the Theorem \ref{efficiencydimensions}.\\
The rest can be analyzed via splitting into the even and odd cases as 
$$x-y=2m, 2m+1  ~~\tn{for}~~ m\geq 1.$$  
We go over two cases for illustrative purposed, their general cases has the same attributes. 
If $x-y=4$ then the coordinate ring can be computed as
$$\mbbc[\{X^4,X^3Y\cdots Y^4\}Z^y, \{X^2,XY,Y^2\}Z^{y+1}W,Z^{y+2}W^2].$$
If $x-y=3$ then the coordinate ring is
$$\mbbc[\{X^3\cdots Y^3\}Z^y, \{X,Y\}Z^{y+1}W].$$
We detect the toric variety from the polytopes of these rings as in the 
Figure \ref{ygeq1}. The general cases are obtained by extending these polytopes accordingly, which clearly does not change the lattice. A straightforward generalization.

\vspace{.5cm}

\end{proof}

\beg{cor}\label{efficiencyrestquotients} 
The quotient is the weighted projective space $\mbb{CP}_{1,1,2}$ for the linearizations in the cases $(2m,0),(2m+1,0)$ for $m\geq 1$ of the Theorem \ref{efficiencydimensions}. \end{cor} 
\beg{proof} The argument of the Theorem \ref{efficiencydimensions2} is still valid in the case of $y=0$.
\end{proof}


In summary, we compute the minimal efficiency dimension to be equal to $1$, 
and this is achieved by the cases $\alpha=(2m,0)$ for $m\geq 1$. In all of these minimal cases we have proved that the quotient is the weighted projective space $\mbb{CP}_{1,1,2}$. We also computed the efficiency and quotients for all the remaining cases.

{\small \beg{flushleft} \textsc{Department of Mathematics, University of Wisconsin at Madison}\\
\textit{E-mail address :} \texttt{\textbf{kalafat@math.wisc.edu}} 


\end{flushleft}

}




  \end{document}